\documentclass[12pt]{article}

\newtheorem{theorem}{Theorem}[section]
\newtheorem{lemma}[theorem]{Lemma}
\newtheorem{proposition}[theorem]{Proposition}
\newtheorem{corollary}[theorem]{Corollary}

\newtheorem{definition}[theorem]{Definition}

\newfam\msbfam
\font\tenmsb=msbm10  scaled \magstep1 \textfont\msbfam=\tenmsb
\font\sevenmsb=msbm7 scaled \magstep1 \scriptfont\msbfam=\sevenmsb
\font\fivemsb=msbm5  scaled \magstep1 \scriptscriptfont\msbfam=\fivemsb
\def\Bbb{\fam\msbfam \tenmsb}

\def\RR{{\Bbb R}}
\def\CC{{\Bbb C}}
\def\QQ{{\Bbb Q}}
\def\NN{{\Bbb N}}
\def\ZZ{{\Bbb Z}}
\def\II{{\Bbb I}}
\def\TT{{\Bbb T}}
\def\BB{{\Bbb B}}

\def\ss{\subseteq}
\def\ra{\rightarrow}

\def\O{\Omega}

\def\Aut{\hbox{Aut}\,}

 \def\HollowBox #1#2{{\dimen0=#1 \advance\dimen0 by -#2       
       \dimen1=#1 \advance\dimen1 by #2                       
        \vrule height #1 depth #2 width #2                    
        \vrule height 0pt depth #2 width #1                   
        \llap{\vrule height #1 depth -\dimen0 width \dimen1}%
       \hskip -#2                                             
       \vrule height #1 depth #2 width #2}}                   
 \def\BoxOpTwo{\mathord{\HollowBox{6pt}{.4pt}}\;}             

\def\endpf{\hfill $\BoxOpTwo$ \smallskip \\ }

\def\dbar{\overline{\partial}}

\newfam\msbbfam
\font\tenmsbb=msbm10  scaled \magstep1 \textfont\msbbfam=\tenmsbb
\font\sevenmsbb=msbm7  scaled \magstep1 \scriptfont\msbbfam=\sevenmsbb
\font\fivemsbb=msbm5    scaled \magstep1 \scriptscriptfont\msbbfam=\fivemsbb

\usepackage{graphicx}

\usepackage{amsmath}

\begin{document}

\begin{center}
\Large \bf On a Theorem of Bers, with Applications to the Study of Automorphism Groups of Domains\footnote{{\bf Subject 
Classification Numbers:}  32A38, 30H50, 32A10, 32M99.}
\end{center}
\vspace*{.12in}

\begin{center}
Steven G. Krantz
\end{center}

\date{\today}

\begin{quote}
{\bf Abstract:}   
We study and generalize a classical theorem of L. Bers that classifies
domains up to biholomorphic equivalence in terms of the algebras of 
holomorphic functions on those domains.  Then we develop applications
of these results to the study of domains with noncompact automorphism groupg.
\end{quote}
\vspace*{.25in}

\markboth{STEVEN G. KRANTZ}{BERS'S THEOREM AND AUTOMORPHISM GROUPS}

\section{Introduction}

For us a {\it domain} in complex space is a connected open set.
If $\Omega$ is a domain then let ${\cal O}(\Omega)$ denote
the algebra of holomorphic functions on $\Omega$.

In 1948, Lipman Bers [BERS] proved the following elegant result:

\begin{theorem} \sl 
Let $\Omega$, $\Omega$ be domains in
$\CC$. If ${\cal O}(\Omega)$ is isomorphic to ${\cal
O}(\widehat{\Omega})$ as an algebra, then the domain $\Omega$ is
conformally equivalent to the domain $\widehat{\Omega}$. 
\end{theorem}

Since that time, this result has been generalized to domains in $\CC^n$, and even
to domains in Stein manifolds---see for instance [ZAM1], [ZAM2].  

In the present paper we offer some other variants of Bers's theorem, and then develop
applications of these results to the study of the automorphism groups of domains
in complex space.

\section{Variants of Bers's Theorem}

In this section we formulate several variants of Bers's theorem.  They
all have the same proof.  For completeness, we provide here the proof
of Bers's original theorem stated in the last section.
\smallskip \\

\noindent {\bf Proof of Theorem 1.1:}  In fact we shall prove the result in $\CC^n$.
Let $\O \ss \CC$ be a domain. Let ${\cal O}(\O)$ denote the algebra of
holomorphic functions from $\O$ to $\CC$. Bers's theorem says, in effect,
that\index{Bers's theorem} 
the algebraic structure of ${\cal O}(\O)$ characterizes $\O$. We
begin our study by introducing a little terminology.

\begin{definition} \rm
Let $\O \ss \CC$ be a domain.  A $\CC$-algebra homomorphism
$\varphi: {\cal O}(\O) \ra \CC$ is called a {\it character}
of ${\cal O}(\O)$.  If $c \in \CC$, then the mapping
\begin{eqnarray*}
e_c: {\cal O}(\O) & \ra & \CC \, , \\
	       f  & \mapsto & f(c) \, , \\
\end{eqnarray*}
is\index{point evaluation} 
called a {\it point evaluation}.  Every point evaluation
is a character.
\end{definition}

It should be noted that, if $\varphi: {\cal O}(\O) \ra {\cal O}(\widehat{\O})$ 
is not the trivial zero homomorphism, then $\varphi(1) = 1$.
This follows because $\varphi(1) = \varphi(1 \cdot 1) = \varphi(1) \cdot \varphi(1)$.
On any open set where the holomorphic function $\varphi(1)$ does not
vanish, we find that $\varphi(1) \equiv 1$.  The result follows
by analytic continuation.

It turns out that every character of ${\cal O}(\O)$ is a point
evaluation.  That is the content of the next lemma.

\begin{lemma} \sl
Let $\varphi$ be a character on ${\cal O}(\O)$.  Then 
$\varphi = e_c$ for some $c \in \O$.  Indeed, $c = \varphi(\hbox{id}) \in \O$.
Here $\hbox{id}$ is defined by $\hbox{id}(z) = z$.
\end{lemma}
\noindent {\bf Proof:}   Let $c$ be defined as in the statement
of the lemma.  Let $f(z) = z - c$.  Then
$$
\varphi(f)  =  \varphi(\hbox{id}) - \varphi(c) 
	    =  c - c 
	    =  0 \, .
$$

If it were not the case that $c \in \O$ then the function
$f$ would be a unit in ${\cal O}(\O)$.  But then
$$
1 = \varphi(f \cdot f^{-1}) = \varphi(f) \cdot \varphi(f^{-1})
     = 0 \, .
$$
That is a contradiction.  So $c \in \O$.

Now let $g \in {\cal O}(\O)$ be arbitrary.  Then we may write
$$
g(z) = g(c) + f(z) \cdot \widetilde{g}(z) \, ,
$$
where $\widetilde{g} \in {\cal O}(\O)$.  Thus
$$
\varphi(g) = \varphi(g(c)) + \varphi(f) \cdot \varphi(\widetilde{g}) 
	   = g(c) + 0 = g(c) = e_c(g) \, .
$$
We conclude that $\varphi = e_c$, as was claimed.
\endpf 
\smallskip \\

Now we may prove Bers's theorem.  We formulate the result in slightly
greater generality than stated heretofore.

\begin{theorem} \sl
Let $\O$, $\widehat{\Omega}$ be domains.  Suppose that
$$
\varphi: {\cal O}(\O) \ra {\cal O}(\widehat{\Omega})
$$
is a $\CC$-algebra homomorphism.  Then there exists
one\index{Bers's theorem} and only one holomorphic mapping $h: \widehat{\Omega} \ra \O$
such that 
$$
\varphi(f) = f \circ h \quad \hbox{for all} \ f \in {\cal O}(\O) \, .
$$
In\index{domain, characterization of in terms of holomorphic function algebra}
fact, the mapping $h$ is given by $h = \varphi(\hbox{id})$.

The homomorphism $\varphi$ is bijective if and only if $h$
is conformal, that is, a one-to-one and onto holomorphic mapping
from $\widehat{\Omega}$ to $\O$.
\end{theorem}
\noindent {\bf Proof:}  Since we want the mapping $h$ to satisfy
$\varphi(f) = f \circ h$ for all $f \in {\cal O}(\O)$, it must
in particular satisfy $\varphi(\hbox{id}_\O) = \hbox{id}_\O \circ h = h$.
We take this as our definition of the mapping $h$.

If $a \in \widehat{\Omega}$, then $e_a \circ \varphi$ is a character
of ${\cal O}(\O)$.  Thus our lemma tells us that
$e_a \circ \varphi$ must in fact be a point evaluation on $\O$.
As a result,
$$
e_a \circ \varphi = e_c \, , \quad \hbox{with} \ 
    c = (e_a \circ \varphi)(\hbox{id}_\Omega) = e_a(h) = h(a) \, .
$$
Thus, if $f \in {\cal O}(\O)$, then
$$
\varphi(f)(a) = e_a(\varphi \circ f) = (e_a \circ \varphi)(f)
     = e_{h(a)}(f) = f(h(a)) = (f \circ h)(a) 
$$
for all $a \in \widehat{\Omega}$.
We conclude that $\varphi(f) = f \circ h$ for all $f \in {\cal O}(\O)$.  

For the last statement of the theorem, suppose that $h$ is a 
one-to-one, onto conformal mapping of $\widehat{\Omega}$ to $\O$.
If $g \in {\cal O}(\O)$, then set $f = g \circ h^{-1}$.  It follows
that $\varphi(f) = f \circ h = g$.  Hence $\varphi$ is onto.
Likewise, if $\varphi(f_1) = \varphi(f_2)$, then $f_1 \circ h = f_2 \circ h$
hence, composing with $h^{-1}$, $f_1 \equiv f_2$.  So $\varphi$ is
one-to-one.  Conversely, suppose that $\varphi$ is an isomorphism.
Let $a \in \O$ be arbitrary.  Then $e_a$ is a character on ${\cal O}(\O)$;
hence $e_a \circ \varphi^{-1}$ is a character on ${\cal O}(\widehat{\Omega})$.
By the lemma, there is a point $c \in \widehat{\Omega}$ such
that $e_a \circ \varphi^{-1} = e_c$.  It follows that
$$
e_a = e_c \circ \varphi \, .
$$
Applying both sides to $\hbox{id}_\O$ yields
$$
e_a(\hbox{id}_\O) = (e_c \circ \varphi)(\hbox{id}_\O) \, .
$$
Unraveling the definitions gives
$$
a = e_c(\hbox{id}_\O \circ h) = h(c) \, .
$$
Thus $h(c) = a$ and $h$ is surjective.  The argument in fact
shows that the pre-image $c$ is uniquely determined.  So $h$
is also one-to-one.
\endpf
\smallskip \\

Now we formulate some variants of Bers's theorem.  Again we stress
that each has the same proof (the proof that we just presented).

In what follows, we shall be dealing with the space $L(\Omega)$ of Lipschitz
functions on $\Omega$.  These are functions that satisfy a condition
of the form
$$
\sup_{x, y \in \Omega \atop x \ne y} \frac{|f(x) - f(y)|}{|x - y|} \leq C \, .  \eqno (2.4)
$$
As usual, we use the expression (2.4) to define a norm $\| \ \ \|_{L(\Omega)}$ on $L(\Omega)$.
	   
\setcounter{theorem}{4}

\begin{proposition} \sl
If $\Omega$ is a  domain in $\CC^n$, then let $L(\Omega)$ denote the
algebra of Lipschitz holomorphic functions on $\Omega$.	 The 
domains $\Omega$ and $\widehat{\Omega}$ in $\CC^n$ are biholomorphically equivalent
if and only if the algebras $L(\Omega)$ and $L(\widehat{\Omega})$ are isomorphic
as algebras.
\end{proposition}

\begin{proposition} \sl 
The bounded domains $\Omega$ and
$\widehat{\Omega}$ in $\CC^n$ are biholomorphically equivalent, with a
biholomorphism that is bi-Lipschitz, if and only if the
algebras $L(\Omega)$ and $L(\widehat{\Omega})$ are isomorphic as
algebras. 
\end{proposition}

We remark that it is possible to formulate versions of these results
for Sobolev spaces of holomorphic functions, for Besov spaces of holomorphic
functions, and in other contexts as well.  We leave the details for the interested
reader.

\section{Applications}

Our intention here is to study the automorphism groups of domains in $\CC^n$.
Here, if $\Omega \ss \CC^n$ is a domain, then the automorphism group
of $\Omega$ (denoted $\hbox{Aut}\, (\Omega)$) is the collection of
biholomorphic mappings of $\Omega$ to itself.  The usual topology
on $\hbox{Aut}\, (\Omega)$ is that of uniform convergence on
compact sets (equivalently, the compact-open topology).   For a bounded
domain $\Omega$, this topology turns $\hbox{Aut}\, (\Omega)$ into
a real Lie group.  Note, however, that the automorphism group
of $\Omega = \CC^n$ with $n > 1$ is infinite dimensional hence certainly {\it not}
a Lie group.

If $\Omega$ is a fixed domain in $\CC^n$ and if $f \in L(\Omega)$, then let us
say that $f$ is {\it noncompact} if there is a sequence
$\varphi_j \in \hbox{Aut}\, (\Omega)$ such that $\{f \circ \varphi_j\}$
is a noncompact set in $L(\Omega)$.  Notice that, obversely, $f$
is compact if $\{f \circ \varphi_j\}$ is a compact set in $L(\Omega)$
for every choice of $\varphi_j$.

\begin{proposition} \sl
Let $\Omega$ be a smoothly bounded, pseudoconvex domain in $\CC^n$.   Then $\Omega$ has noncompact
automorphism group if and only if there exists an $f \in L(\Omega)$ such that $f$
is noncompact.
\end{proposition}
{\bf Proof:}  If the automorphism group is noncompact, then (by a classical
result of H. Cartan), there exist $\varphi_j \in \hbox{Aut}\, (\Omega)$,
$P \in \Omega$, and $X \in \partial \Omega$ such that $\varphi_j(P) \ra X$.
By a result of Ohsawa (see [OHS]), the Bergman metric is complete.
Fix a nonconstant $f \in L(\Omega)$.  Choose $p, q \in \Omega$, $p\ne q$, so that
$$
|p - q| \approx (1/|f\|_{L(\Omega)}) \cdot |f(p) - f(q)| \, .
$$
We may suppose without loss of generality that $|p - q| = 1$. 

Now certainly $|\varphi_j(p) - \varphi_j(q)| \ra 0$
(since, by the completeness of the metric, both $\varphi_j(p)$ and $\varphi_j(q)$ must both tend to $X$).  
We may now calculate that
\begin{eqnarray*}
       C  & = & C |p - q|  \\
        & \approx & (1/\|f\|_{L(\Omega)}) \cdot |f(p) - f(q)|  \\
         & = & (1/\|f\|_{L(\Omega)}) \cdot |f(\varphi^{-1}_j(\varphi_j(p))) - f(\varphi^{-1}_j(\varphi_j(q)))| \, .
\end{eqnarray*}
Since $|\varphi_j(p) - \varphi_j(q)| \ra 0$, we see that $\{f \circ \varphi^{-1}_j\}$ has
Lipschitz norm which is blowing up.  So $f$ is noncompact.

Conversely, if $\Aut(\Omega)$ is compact, then let $f \in L(\Omega)$ and consider
$\{f \circ \varphi_j\}$ for $\varphi_j \in \Aut(\Omega)$.  Examine
$$
|f \circ \varphi_j(p) - f \circ \varphi_j(q)| \, .   \eqno (3.1.1)
$$
Clearly, by compactness, $|\nabla \varphi_j|$ is bounded above and below,
uniformly in $j$, on any compact set $K \subset \! \subset \Omega$.   By the Ascoli-Arzela theorem applied on compact sets, we
see from (3.1.1) that $f \circ \varphi_j$ has a convergent subsequence.
\endpf 
\smallskip \\

The next well-known result, due to Bun Wong [WON], is a cornerstone
of the modern theory of automorphism groups of smoothly bounded domains.
Now we present some new proofs of this result.

\begin{theorem} \sl
Let $\Omega$ be a smoothly bounded, strongly pseudoconvex domain in $\CC^n$.  Suppose that
there a point $P \in \Omega$ and a strongly pseudoconvex
boundary point $X \in \partial \Omega$ and that
there exist $\varphi_j \in \Aut (\Omega)$ such that
$\varphi_j(P) \ra X$.  Then $\Omega$ is biholomorphic to the unit
ball $B \ss \CC^n$.
\end{theorem}
{\bf Proof:}  As advertised, we shall sketch three proofs.  We first
note that, according to H. Cartan's theorem and our previous result,
the hypotheses imply that there is an $f \in L(\Omega)$ which is noncompact.
\medskip \\

\noindent {\bf First Proof of the Theorem:}  If $\Omega$ is {\it not}
biholomorphic to the ball then, by a celebrated result of Lu Qi-Keng [LQK]
(see [GKK] for thorough discussion), there is a point $Q$ in $\Omega$
where the holomorphic sectional curvature of the Bergman metric
is not the constant holomorphic sectional curvature of the ball.

As noted in the proof of the preceding result, the Bergman metric
is complete on $\Omega$.   So in fact any compact set $K \subset \! \subset \Omega$ has the 
property that $\{\varphi_j\}$ converges uniformly on $K$ to $X$.
In particular, $\varphi_j(Q) \ra X$.  But it can be calculated
(see [KLE], [GK1], [GKK]) that the holomorphic sectional
curvature of the Bergman metric tends to the constant curvature
of the ball at points that approach a strongly pseudoconvex boundary
point $X$.   That contradicts the last sentence of the previous paragraph.

We conclude that $\Omega$ {\it is} biholomorphic to the ball, as claimed.
\endpf
\medskip \\

\noindent {\bf Second Proof of the Theorem:}  
It is convenient for this argument to equip ${\cal O}(\Omega)$ with the topology
of uniform convergence on compact sets (i.e., the compact-open topology).  For
convenience, and without any loss of generality, we restrict attention now
to ambient dimension 2.

Let $U$ be a small neighborhood of $X$.  Since $X$ is a peak point (see [KRA1]), it is standard to argue that,
for any compact set $K \ss \Omega$, there is a $J$ so large that $j > J$ implies
that $\varphi_j(K) \ss U \cap \Omega$.  Let $X'$ be a point of $U \cap \Omega$ that
is very near to $X$.  Let $\delta = \delta_j = \hbox{dist}\,(X', \partial \Omega)$.  
After a normalization of coordinates, we may suppose that the complex normal direction
at $X$ is $z_1$ and the complex tangential direction at $X$ is $z_2$.

Define $\psi(z_1, z_2) = (X'_1 + (z_1 - X'_1)/\delta, X'_2 + (z_2 - X'_2)/\sqrt{\delta})$.
Then $\psi \circ \varphi_j$, with $j$ as above, will have Lipschitz norm that is bounded, 
independent of $j$.  As a result, using a sequence of compact sets $K_j$ that exhausts $\Omega$,
and neighborhoods $U$ that shrink to $X$,
we may derive a subsequence, convergent on compact sets.
And it will converge to a mapping of $\Omega$ to the Siegel upper half space.  [This is just
the standard method of scaling, which is described in detail in [GKK]).  So $\Omega$ is
biholomorphic to the Siegel upper half space, which is in turn biholomorphic
to the unit ball.
\endpf
\medskip \\

\noindent {\bf Third Proof of the Theorem:}  For this proof we examine
the Fefferman asymptotic expansion for the Bergman kernel near a 
strongly pseudoconvex boundary point (see [FEF] and also [GKK]).
This says that, in suitable local coordinates, 
$$
K(z, \zeta) = \frac{\psi(z, \zeta)}{[- X(z, \zeta)]^{n+1}} + \widetilde{\psi}(z, \zeta) \cdot \log [- X](z, \zeta) \, .  \eqno (3.2.1)
$$
Here $\psi$, $\widetilde{\psi}$ are smooth functions on $\overline{\Omega} \times \overline{\Omega}$ and
$X$ is the Levi polynomial (see [KRA1, Ch.\ 3]) on $\Omega$.

An interesting feature of Fefferman's work, and subsequent work of Burns and Graham [GRA], is
that the logarithmic term is always present near a boundary point that is not spherical.

Arguing as usual, if $P$ and $X$ exist then any other point $Q
\in \Omega$ has the property that $\varphi_j(Q) \ra X$ as $j
\ra \infty$. We begin with a point $Q$ near the boundary at
which the Fefferman expansion (3.2.1) is valid. If $\Omega$
is not the ball then we can take $Q$ to be very near to a
boundary point that is not spherical.

Of course the Bergman kernel transforms under a biholomorphic mapping $F$
of $\Omega$ by the standard formula ([KRA1, Ch.\ 1])
$$
\hbox{Jac}_\CC F(z) K(F(z), F(\zeta)) \overline{\hbox{Jac}_\CC F(\zeta)} = K(z, \zeta) \, .  \eqno (3.2.2)
$$
So, when we think of $\varphi_j(Q) \ra X$, then we may understand how
the Bergman kernel transforms by applying the transformation formula (3.2.2)
to the Fefferman expansion (3.2.1).  On the one hand, this should
give rise to another Fefferman-type formula based at the point $\varphi_j(Q)$.  But
the problem is that the logarithmic expression does not scale.  The result, as $j \ra \infty$,
will not be a valid Fefferman formula.  That is a contradiction.  So $\Omega$ must
be biholomorphic to the ball.
\endpf
\smallskip \\

A consequence of the first two results is this:

\begin{corollary} \sl
A strongly pseudoconvex domain $\Omega \ss \CC^N$ is biholomorphic
to the ball if and only if the algebra $L(\Omega)$ of Lipschitz 
functions is noncompact.
\end{corollary}

\section{An Analysis of Algebra Isomorphisms}

First suppose that $A$ is an annulus in the complex plane.  Suppose
that
$$
\Phi: {\cal O}(A) \longrightarrow {\cal O}(A) 
$$
is an algebra isomorphism.  We claim that $\Phi(z) = z$.  That is to say,
$\Phi$ maps the holomorphic identity function to itself.

First of all, it cannot be that $\Phi(z) = z^2$ or any other higher-order
polynomial (or power series or Laurent series) because then it is clear
that $\Phi$ would not be onto.	A similar argument shows that $\Phi(z)$ cannot
be a Laurent series with initial term having negative index.

So $\Phi(z)$ is either a power series beginning with a degree-zero term
or a power series beginning with a degree-one term.  But
obviously $\Phi(1) = 1$ and $\Phi(0) = 0$. So $\Phi(z)$ is a power series
beginning with a first-degree term.   But in fact if that power series contains
any term beyond the first-degree term, then there is no holomorphic function that
will map to $z^2$ under $\Phi$.  So the power series is simply of the form $\alpha z$.
So any Laurent series of the form
$$
\sum_{j = -\infty}^\infty a_j z^j  \eqno (4.1)
$$
is mapped under $\Phi$ to
$$
\sum_{j=-\infty}^\infty \alpha^j a_j z^j \, .  \eqno (4.2)
$$
If the Laurent series in (4.1) is chosen so that the function
it defines has $A$ as its natural domain of definition, and if the modulus
of $\alpha$ is not 1, then it follows that the image function given
by (4.2) will have a {\it different} natural domain of
definition.  And that is impossible.

We conclude that $\alpha$ has modulus 1.  We may as well take $\alpha = 1$.

This example illustrates Theorem 3.2.  For the automorphism group of an
annulus is just two copies of the circle group.  So it is compact.  As a result,
any $f \in L(A)$ will be compact.

\section{Further Results}

The next result is classical.	See [KRA2, Ch.\ 12] for a more traditional proof.

\begin{proposition} \sl
Fix a bounded domain $\Omega \ss \CC^n$.  Let $\{\varphi_j\}$ be automorphisms
of $\Omega$.  Assume that the $\varphi_j$ converge normally (i.e., uniformly on
compact sets) to a limit $f$.  Then either
\begin{enumerate}
\item[{\bf (1)}]  The mapping $f$ is an automorphism of $\Omega$;

or

\item[{\bf (2)}]  The mapping $f$ is a constant.
\end{enumerate}
\end{proposition}
{\bf Proof:}  We adopt the point of view of Bers's theorem.

With $\varphi_j \in \Aut(\Omega)$ as in the statement of the theorem, and
$g \in L(\Omega)$, examine $\{g \circ \varphi_j\}$.   

Now either $g \circ \varphi_j$ is compact or it is not. If $g
\circ \varphi_j$ is compact, then there exists a subsequence $\varphi_{j_k}$ and a $\tau$ such that
$g \circ \varphi_{j_k} \ra \tau$ with $\tau \in L(\Omega)$. So $g \circ f = \tau$, with $f \in \Aut(\Omega)$ (because
it is a nondegenerate mapping, and a limit of automorphisms).  Specifically, the mapping $f$ is univalent because 
it is the limit of univalent mappings.  Also $f$ is onto because we can apply 
our reasoning to $\varphi_j^{-1}$.  That
is part {\bf (1)} of our conclusion (formulated in the language of the present paper).

If instead $g \circ \varphi_j$ is noncompact, then $\{g \circ \varphi_j\}$ has no convergent
subsequence.  So $g \circ \varphi_j$ blows up in norm.  Hence there are a point $P \in \Omega$
and a point $X \in \partial \Omega$ such that $\varphi_{j_k}(P) \ra X$ (for some subsequence
$\varphi_{j_k}$).  Hence $g \circ \varphi_{j_k} \ra g(X)$.  That completes the proof of {\bf (2)}.
\endpf
\smallskip \\

Now we have

\begin{proposition} \sl
Suppose that $f:\Omega \ra \Omega$ is a holomorphic mapping.  Assume that,
for some sequence $\{\varphi_j\}$ of automorphisms of $\Omega$, $f \circ \varphi_{j_k}$
converges normally to a function $g \in {\cal O}(\Omega)$.  Then
\begin{enumerate}
\item[{\bf (a)}]  If $g \in \Aut(\Omega)$, then $f \in \Aut(\Omega)$;
\item[{\bf (b)}]  If $g$ is not constant then every convergent subsequence
of $h_k \equiv f \circ \varphi_{j_{k+1}} \circ \varphi_{j_k}^{-1}$ has limit $\hbox{id}_\Omega$.
\end{enumerate}
\end{proposition}

{\bf Proof:}  This result is like a converse to compactness.

If $f(a) = f(b)$ for some distinct points $a, b \in \Omega$ then
$$
f(\varphi_{j_k} \circ \varphi_{j_k}^{-1}(a)) = f(\varphi_{j_k} \circ \varphi_{j_k}^{-1}(b))  \, .
$$
Now, if the $\varphi_{j_k}$ converge to some $\psi$, then we see that
$$
g(\psi(a)) = g(\psi(b)) \, .
$$
If $\psi$ is an automorphism then this is certainly a contradiction.

Of course $f \circ \varphi_{j_k}(\Omega) \ss f(\Omega)$ for all $k$.
So $g(\Omega) \ss f(\Omega) \ss \Omega$.  But $g(\Omega) = \Omega$.  So $f(\Omega) = \Omega$.
Thus $f$ is onto.  It is also one-to-one.  This proves {\bf (a)}.

For part {\bf (b)}, we take $g$ to be holomorphic and nonconstant.  Let $h$ be a subsequential
limit of $f \circ \varphi_{j_{k+1}} \circ \varphi_{j_k}^{-1} \equiv h_k$.
As a result, $f \circ \varphi_{j_{k+1}} = h_k \circ \varphi_{j_k}$ so $g = h \circ \psi$.
But then $h = g \circ \psi^{-1}$.  So $h$ differs from $g$ by an automorphism.  Certianly
then $h$ is nonconstant.  We note further that $g = h \circ \psi$ so that $g$ is
an automorphism.
\endpf
\smallskip \\

\newpage

\noindent {\Large \sc References}
\vspace*{.2in}

\newfam\msbfam
\font\tenmsb=msbm10  scaled \magstep1 \textfont\msbfam=\tenmsb
\font\sevenmsb=msbm7 scaled \magstep1 \scriptfont\msbfam=\sevenmsb
\font\fivemsb=msbm5  scaled \magstep1 \scriptscriptfont\msbfam=\fivemsb
\def\Bbb{\fam\msbfam \tenmsb}

\def\RR{{\Bbb R}}
\def\CC{{\Bbb C}}
\def\QQ{{\Bbb Q}}
\def\NN{{\Bbb N}}
\def\ZZ{{\Bbb Z}}
\def\II{{\Bbb I}}
\def\TT{{\Bbb T}}
\def\BB{{\Bbb B}}

\begin{enumerate}

\item[{\bf [BERS]}] L. Bers, On rings of analytic functions,
{\it Bull.\ Amer.\ Math.\ Soc.} 54(1948), 311--315.

\item[{\bf [FEF]}] C. Fefferman, The Bergman kernel and
biholomorphic mappings of pseudoconvex domains, {\it Invent.
Math.} 26(1974), 1--65.
			     
\item[{\bf [GRA]}] C. R. Graham, Scalar boundary invariants and
the Bergman kernel, {\it Complex analysis}, II (College Park,
Md., 1985--86), 108--135, {\it Lecture Notes in Math.}, 1276,
Springer, Berlin, 1987.

\item[{\bf [GKK]}]  R. E. Greene, K.-T. Kim, and S. G. Krantz,
{\it The Geometry of Complex Domains},
Birkh\"{a}user Publishing, Boston, 2011.

\item[{\bf [GK1]}]  R. E. Greene and S. G. Krantz, Deformations of complex structure, 
estimates for the  $\dbar$-equation, and stability of the Bergman
kernel, {\it Advances in Math.} 43(1982), 1--86.
					      
\item[{\bf [KLE]}]  P. Klembeck, K\"{a}hler metrics of negative curvature, the
Bergman metric near the boundary and the Kobayashi metric on
smooth bounded strictly pseudoconvex sets, {\it Indiana Univ.\ Math.\ J.}
27(1978), 275--282. 

\item[{\bf [KRA1]}]  S. G. Krantz, {\it Function Theory of
Several Complex Variables}, 2nd. ed., American Mathematical
Society, Providence, RI, 2001.

\item[{\bf [KRA2]}] S. G. Krantz, {\it Cornerstones of
Geometric Function Theory: Explorations in Complex Analysis},
Birkh\"{a}user Publishing, Boston, 2006.

\item[{\bf [LQK]}] Lu Qi-Keng, On K\"{a}hler manifolds with
constant curvature, {\it Acta.\ Math.\ Sinica} 16(1966),
269--281 (Chinese) (= {\it Chinese Math.} 9(1966), 283--298).

\item[{\bf [OHS]}] T. Ohsawa, A remark on the completeness of
the Bergman metric, {\it Proc.\ Japan Acad.} Ser.\ A Math.
Sci. 57(1981), 238--240.

\item[{\bf [WON]}] B. Wong, Characterizations of the ball in
$\CC^n$ by its automorphism group, {\it Invent.\ Math.}
41(1977), 253--257.

\item[{\bf [ZAM1]}] W. R. Zame, Homomorphisms of rings of germs
of analytic functions, {\it Proc.\ Amer.\ Math.\ Soc.}
33(1972), 410--414.

\item[{\bf [ZAM2]}] W. R. Zame, Induced homomorphisms of
algebras of analytic germs, {\it Complex Analysis}, 1972
(Proc.\ Conf., Rice Univ., Houston, Tex., 1972), Vol. II:
Analysis on singularities. Rice Univ.\ Studies 59 (1973), no.
2, 157--163.

\end{enumerate}
\vspace*{.25in}

\begin{quote}
Steven G. Krantz  \\
Department of Mathematics \\
Washington University in St. Louis \\
St.\ Louis, Missouri 63130  \\
{\tt sk@math.wustl.edu}
\end{quote}

\end{document}